\documentclass[12pt]{amsart}
\usepackage[margin=1.2in]{geometry}
\usepackage[all,cmtip]{xy}

\theoremstyle{plain}
\newtheorem{theorem}{Theorem}
\newtheorem{prop}[theorem]{Proposition}
\newtheorem{corollary}[theorem]{Corollary}
\newtheorem{lemma}[theorem]{Lemma}
\newtheorem{claim}[theorem]{Claim}
\newtheorem{conj}[theorem]{Conjecture}

\theoremstyle{definition}
\newtheorem{defn}[theorem]{Definition}
\newtheorem{rmk}[theorem]{Remark}

\numberwithin{theorem}{section}

\newcommand{\spec}{\operatorname{Spec}}
\newcommand{\Hom}{\operatorname{Hom}}
\newcommand{\ord}{\operatorname{ord}}

\begin{document}

\title{One higher dimensional analog of jet schemes}

\author[C. Yuen]{Cornelia Yuen} 
\address{SUNY Potsdam, Department of Mathematics, 44 Pierrepont
Avenue, Potsdam, NY 13676, USA} 
\email{{\tt yuenco@potsdam.edu}}

\begin{abstract}
    We develop the theory of truncated wedge schemes, a higher
    dimensional analog of jet schemes.  We prove some basic properties
    and give an irreducibility criterion for truncated wedge schemes
    of a locally complete intersection variety analogous to
    Musta\c{t}\v{a}'s for jet schemes.  We show that the reduced
    subscheme structure of a truncated wedge scheme of any monomial
    scheme is itself a monomial scheme in a natural way by explicitly
    describing generators of each of the minimal primes of any
    truncated wedge scheme of a monomial hypersurface.  We give
    evidence that the irreducible components of the truncated wedge
    schemes of a reduced monomial hypersurface all have multiplicity
    one.
\end{abstract}

\thanks{This is part of the author's PhD thesis at the University of
Michigan under the supervision of Karen Smith.  The author was
supported by her advisor's NSF grant $0070722$ and the department's
NSF RTG grant $0502170$.}

\maketitle

\section{Introduction}

A wedge is a two-dimensional analog of an arc. More precisely, for 
$X$ a scheme of finite type over an algebraically closed field $k$ of 
characteristic zero, a wedge on $X$ is a $k$-morphism $\spec k[[s,t]] 
\rightarrow X$, which can be thought of as an ``infinitesimal surface 
germ'' on $X$.

The study of wedges was initiated by Lejeune-Jalabert in $1980$ in an
attack on Nash's conjecture \cite{Lejeune-Jalabert}.  Her idea rested
on the observation that any wedge $\spec k[[s,t]] \rightarrow X$ can
be precomposed with the natural map $\spec k[[t]] \rightarrow \spec
k[[s,t]]$, dual to the map of rings sending $t \mapsto t$ and $s
\mapsto t$, producing an arc called the ``center'' of the wedge.  She
showed that Nash's conjecture (for surfaces) could be settled by an
affirmative answer to the following question: Does a wedge centered at
a ``general'' arc on a normal surface singularity lift to its minimal
resolution of singularities?  Later, Reguera \cite{RegueraCURVE}
considered this problem for wedges on higher dimensional varieties,
and showed that a positive answer to Nash's question is equivalent to
a positive answer to this wedge extension problem.

In this paper, we study the analogous notion to $m$-jets. Recall that 
an $m$-jet is a $k$-morphism $\spec k[t]/(t^{m+1}) \rightarrow X$. It 
is known that the set of all $m$-jets of $X$ carries a natural scheme 
structure, called the $m^{th}$ jet scheme of $X$ and denoted by 
$\mathcal{J}_{m}(X)$\footnote{For a discussion of the fundamentals of 
jet schemes, we refer the readers to \cite{Blickle}.}.

\begin{defn}
    Let $X$ be a scheme of finite type over $k$. An $m$-wedge of $X$ 
    is a $k$-morphism $$\spec k[s,t]/(s,t)^{m} \rightarrow X.$$
\end{defn}

As with $m$-jets, the collection of all $m$-wedges forms a scheme
$\mathcal{W}_{m}(X)$ in a natural way, called the $m^{th}$ wedge
scheme of $X$.  We also have the natural projection maps
$\pi_{m-1}^{m}:\mathcal{W}_{m}(X) \rightarrow \mathcal{W}_{m-1}(X)$
and $\pi_{m}:\mathcal{W}_{m}(X) \rightarrow X$ induced by pulling back
an $m$-wedge via the natural maps $\spec k[s,t]/(s,t)^{m+1}
\hookrightarrow \spec k[s,t]/(s,t)^{m}$ and $\spec k \hookrightarrow
\spec k[s,t]/(s,t)^{m+1}$.

Analogous to the situation with arcs, the set of all wedges also forms
a scheme
$\mathcal{W}_{\infty}(X)=\underleftarrow{\lim}\mathcal{W}_{m}(X)$,
called the wedge scheme of $X$.

We begin our study of truncated wedge schemes by showing some basic
properties, including a functorial representation, base change under
\'etale morphisms, smoothness of truncated wedge schemes of a smooth
scheme, and a description of the first truncated wedge scheme in terms
of the first jet scheme.  We also give an irreducibility criterion for
truncated wedge schemes of a locally complete intersection variety
analogous to Musta\c{t}\v{a}'s for jet schemes.  

Next we conduct a detailed study of the truncated wedge schemes of
monomial schemes.  We show that the reduced subscheme structure of a
truncated wedge scheme of any monomial scheme is itself a monomial
scheme in a natural way.  We give explicit generators of each of the
minimal primes of any truncated wedge scheme of a monomial
hypersurface, analogous to that of Goward and Smith for jet schemes.
Lastly, we give evidence that the irreducible components of the
truncated wedge schemes of a reduced monomial hypersurface all have
multiplicity one.

\section{Properties}

First we discuss how the set of all $m$-wedges of $X$ forms a scheme.

Let $X$ be the affine space $\mathbb{A}^{r}$, then an $m$-wedge of $X$
corresponds to a $k$-algebra homomorphism
\begin{align*}
    k[x_{1},\ldots,x_{r}] &\rightarrow k[s,t]/(s,t)^{m+1}\\
    x_{i} &\mapsto x_{i}^{(0,0)}+\ldots
    +x_{i}^{(m,0)}s^{m}+x_{i}^{(m-1,1)}s^{m-1}t+\ldots
    +x_{i}^{(0,m)}t^{m}.
\end{align*}
So the truncated wedge scheme $\mathcal{W}_{m}(\mathbb{A}^{r}) \cong
\spec k[x_{k}^{(i_{k},j_{k})}]$ where $1 \leq k \leq r$ and $0 \leq
i_{k}+j_{k} \leq m$.  In other words, the scheme
$\mathcal{W}_{m}(\mathbb{A}^{r})=\mathbb{A}^{\frac{1}{2}r(m+1)(m+2)}$.

Now if $X \subseteq \mathbb{A}^{r}$ is a closed subscheme, then its
truncated wedge scheme $\mathcal{W}_{m}(X)$ is a closed subscheme of
$\mathcal{W}_{m}(\mathbb{A}^{r})$. Say $X=\spec
k[x_{1},\ldots,x_{r}]/(f_{1},\ldots,f_{d})$. Then an
$m$-wedge of $X$ corresponds to a $k$-algebra homomorphism
\begin{align*}
    \phi:k[x_{1},\ldots,x_{r}]/(f_{1},\ldots,f_{d}) &\rightarrow
    k[s,t]/(s,t)^{m+1}\\
    x_{i} &\mapsto x_{i}^{(0,0)}+\ldots +x_{i}^{(m,0)}s^{m}+\ldots
    +x_{i}^{(0,m)}t^{m},
\end{align*}
subject to the relations $\phi(f_{k})=0$.  Therefore,
$\mathcal{W}_{m}(X)$ is defined by the ideal $W_{m}(X)$ whose
generators $g_{ij}$ are the coefficients of $s^{i}t^{j}$ in
$\phi(f_{k})$ for $1 \leq k \leq d$, $0 \leq i+j \leq m$.  Note that
this computation commutes with localization.  So this local
construction of truncated wedge schemes can be patched together to
give a scheme structure on the set of $m$-wedges of $X$ for
any scheme $X$ of finite type over $k$.

Next, we show that truncated wedge schemes have nice properties
similar to those of jet schemes.  The first example is truncated wedge
schemes as representing schemes of a functor.

\begin{prop}\label{WedgeFunctor}
    Given a scheme $X$ of finite type over $k$, the $m^{th}$ wedge 
    scheme $\mathcal{W}_{m}(X)$ represents the functor
    \begin{align*}
	F:{\bf k}\text{\bf-Schemes} &\rightarrow \text{\bf Sets}\\
	Z &\mapsto \Hom_{k}(Z \times_{k} \spec k[s,t]/(s,t)^{m+1},X),
    \end{align*}
    where the right hand side denotes the set of morphisms of 
    $k$-schemes from the scheme $Z \times_{k} \spec 
    k[s,t]/(s,t)^{m+1}$ to $X$.
\end{prop}
\begin{proof}
    Our goal is to show that the functor $F$ described in the
    proposition is equivalent to the functor of points of
    $\mathcal{W}_{m}(X)$:
    \begin{align*}
	F_{\mathcal{W}_{m}(X)}:{\bf k}\text{\bf-Schemes} &\rightarrow
	\text{\bf Sets}\\
	Z &\mapsto \Hom_{k}(Z,\mathcal{W}_{m}(X)).
    \end{align*}
    
    To understand the set $\Hom_{k}(Z \times_{k} \spec
    k[s,t]/(s,t)^{m+1},X)$, we may assume both $X$ and $Z$ are affine
    \cite[Proposition VI-2]{EisenbudHarris}.  Say $X=\spec
    k[x_{1},\ldots,x_{r}]/(f_{1},\ldots,f_{d})$ and $Z=\spec R$ for
    some $k$-algebra $R$.  An element of $\Hom_{k}(Z \times_{k} \spec
    k[s,t]/(s,t)^{m+1},X)$ corresponds to a $k$-algebra map
    \begin{align*}
	\phi:k[x_{1},\ldots,x_{r}]/(f_{1},\ldots,f_{d}) &\rightarrow
	R[s,t]/(s,t)^{m+1}\\
	x_{k} &\mapsto x_{k}^{(0,0)}+\ldots
	+x_{k}^{(m,0)}s^{m}+x_{k}^{(m-1,1)}s^{m-1}t+\ldots
	+x_{k}^{(0,m)}t^{m}
    \end{align*}
    where $x_{k}^{(i_{k},j_{k})} \in R$ arbitrary, subject to the
    conditions $\phi(f_{l})=0$ for $l=1,\ldots,d$.  So an element of
    $\Hom_{k}(Z \times_{k} \spec k[s,t]/(s,t)^{m+1},X)$ corresponds to
    an $N$-tuple
    $\alpha=(\alpha_{1}^{(0,0)},\ldots,\alpha_{1}^{(1,m-1)},
    \alpha_{1}^{(0,m)},\ldots,\alpha_{r}^{(0,0)},\ldots,
    \alpha_{r}^{(0,m)}) \in R^{N}$ satisfying the equations
    $\phi(f_{l})=0$ for all $l$.  This means that $\alpha$ is an
    $R$-valued point of $\mathcal{W}_{m}(X)$.  In other words, we have
    shown that $$\Hom_{k}(Z \times_{k} \spec
    k[s,t]/(s,t)^{m+1},X)=\Hom_{k}(Z,\mathcal{W}_{m}(X))$$ for all
    $k$-schemes $Z$, and hence the scheme $\mathcal{W}_{m}(X)$
    represents the functor $F$.
\end{proof}

Now using this functorial point of view, we will show that truncated
wedge schemes behave well under \'etale morphisms.
    
\begin{prop}\label{WedgeEtale}
    If $f:X \rightarrow Y$ is an \'etale morphism of finite type
    $k$-schemes, then $\mathcal{W}_{m}(X) \cong \mathcal{W}_{m}(Y)
    \times_{Y} X$ for all $m$.
\end{prop}
\begin{proof}
    We show the equality on the level of the corresponding functor of
    points.  That is, we show the two sets of $k$-morphisms
    $$\Hom_{k}(Z,\mathcal{W}_{m}(X))=\Hom_{k}(Z \times_{k} \spec
    k[s,t]/(s,t)^{m+1},X)$$ and
    \begin{align*}
	&\Hom_{k}(Z,\mathcal{W}_{m}(Y) \times_{Y} X)\\
	= &\Hom_{k}(Z,\mathcal{W}_{m}(Y)) \times_{\Hom_{k}(Z,Y)}
	\Hom_{k}(Z,X)\\
	= &\Hom_{k}(Z \times_{k} \spec k[s,t]/(s,t)^{m+1},Y)
	\times_{\Hom_{k}(Z,Y)} \Hom_{k}(Z,X)
    \end{align*}
    are the same for all $k$-schemes $Z$.  Note it suffices to
    check the equality for all affine schemes $Z$.  Fix a $k$-scheme
    $Z=\spec R$ and consider the commutative diagram
    $$\xymatrix{Z\ar[d]_{\alpha}\ar[r] & Z \times_{k} \spec
    k[s,t]/(s,t)^{m+1}\ar@{-->}[dl]_{\gamma}\ar[d]^{\beta}\\
      X\ar[r]^{f} & Y}$$
    where the top horizontal map is the closed embedding $\spec R
    \hookrightarrow \spec R[s,t]/(s,t)^{m+1}$.  
    
    Obviously, a $Z$-valued $m$-wedge $\gamma \in \Hom_{k}(\spec
    R[s,t]/(s,t)^{m+1},X)$ of $X$ induces a $Z$-valued $m$-wedge
    $\beta \in \Hom_{k}(\spec R[s,t]/(s,t)^{m+1},Y)$ and a map $\alpha
    \in \Hom_{k}(\spec R,X)$ by composition.  Conversely, the maps
    $\alpha$ and $\beta$ together induce a unique map $\gamma$ by the
    definition of formal \'etaleness \cite[Definition (17.1.1)]{EGA}.
    So we have shown
    \begin{align*}
	&\Hom_{k}(Z \times_{k} \spec k[s,t]/(s,t)^{m+1},X)\\
	=&\Hom_{k}(Z \times_{k} \spec k[s,t]/(s,t)^{m+1},Y)
	\times_{\Hom_{k}(Z,Y)} \Hom_{k}(Z,X),
    \end{align*}
    or equivalently,
    $$\Hom_{k}(Z,\mathcal{W}_{m}(X))=\Hom_{k}(Z,\mathcal{W}_{m}(Y)
    \times_{Y} X)$$ for all $k$-schemes $Z$.  That is, the schemes
    $\mathcal{W}_{m}(X)$ and $\mathcal{W}_{m}(Y \times_{Y} X)$ are
    isomorphic.
\end{proof}

When $X$ is a smooth scheme over $k$, we also have the result that 
$\mathcal{W}_{m}(X)$ is an affine bundle over $X$:

\begin{corollary}\label{SmoothWedge}
    Let $X$ be a smooth scheme over $k$ of dimension $n$.  Then
    $\mathcal{W}_{m}(X)$ is locally an
    $\mathbb{A}^{\frac{1}{2}nm(m+3)}$-bundle over $X$.  In particular,
    $\mathcal{W}_{m}(X)$ is smooth of dimension
    $\frac{1}{2}n(m+1)(m+2)$.
\end{corollary}
\begin{proof}
    Since $X$ is smooth over $k$, $X$ is covered by an open affine
    cover $\{U_{i}\}_{i}$ with $U_{i} \rightarrow V_{i}$ \'etale, for
    some $V_{i}$ open subset of $\mathbb{A}^{n}$.  Then
    \begin{align*}
	\mathcal{W}_{m}(U_{i}) &\cong \mathcal{W}_{m}(V_{i})
	\times_{V_{i}} U_{i} \hspace{0.5cm}\text{(by Proposition
	\ref{WedgeEtale})}\\
	&\cong (\mathcal{W}_{m}(\mathbb{A}^{n})
	\times_{\mathbb{A}^{n}} V_{i}) \times_{V_{i}} U_{i}
	\hspace{0.5cm}\text{(since an open immersion is \'etale)}\\
	&\cong \mathcal{W}_{m}(\mathbb{A}^{n}) \times_{\mathbb{A}^{n}}
	U_{i}.
    \end{align*}
    Now $\mathcal{W}_{m}(\mathbb{A}^{n})$ is an
    $\mathbb{A}^{\frac{1}{2}nm(m+3)}$-bundle over $\mathbb{A}^{n}$,
    therefore $\mathcal{W}_{m}(X)$ is an
    $\mathbb{A}^{\frac{1}{2}nm(m+3)}$-bundle over $X$.
\end{proof}

Since a tangent vector on a scheme $X$ is simply a $1$-jet by
definition, the first jet scheme of $X$ is the total tangent space of
$X$.  Interestingly, the first wedge scheme is also related
to it.

\begin{prop}\label{W_{1}(X)}
    Let $X$ be any scheme of finite type over $k$.  Then
    $\mathcal{W}_{1}(X)$ is isomorphic to $\mathcal{J}_{1}(X) \times_{X}
    \mathcal{J}_{1}(X)$ in the category of $k$-schemes.
\end{prop}
\begin{proof}
    We show equality on the level of functors of points; that is, we 
    show the two sets of $k$-morphisms 
    $$\Hom_{k}(Z,\mathcal{W}_{1}(X))=\Hom_{k}(Z \times_{k} \spec 
    k[s,t]/(s,t)^{2},X)$$ and $$\Hom_{k}(Z,\mathcal{J}_{1}(X) \times_{X} 
    \mathcal{J}_{1}(X))$$ are the same for all $k$-schemes $Z$. 
    Again, it is sufficient to check the equality for all affine 
    schemes $Z$.
    
    In the category of {\bf k}\text{\bf-Algebras}, we always have
    coproducts and pushout squares while products and pullback squares
    are much less common (for definition of pushout squares and
    pullback squares, see \cite[p.  65 and p.  71
    respectively]{MacLane}).  However, we do have the following
    pullback square:
    $$\xymatrix{k[s,t]/(s,t)^{2} \ar[d]\ar[r] &
      k[s,t]/(s^{2},t)\ar[d]\\
      k[s,t]/(s,t^{2})\ar[r] & k}$$
    where the horizontal maps are natural surjections sending $t$ to
    zero, and the vertical maps are natural surjections sending $s$ to
    zero.  That is, the ring $k[s,t]/(s,t)^{2}$ is a product in the
    category of $k$-algebras.
    
    We leave the checking to the reader, but caution that analogous
    diagrams with $m>2$ in place of $2$ are {\it not} also pullback
    squares.  Now we apply the contravariant functor $\spec (-)$ to
    our pullback square of $k$-algebras, and get the commutative
    diagram:
    $$\xymatrix{\spec k[s,t]/(s,t)^{2} & \spec
      k[s,t]/(s^{2},t)\ar[l]\\
      \spec k[s,t]/(s,t^{2})\ar[u] & \spec k\ar[l]\ar[u]}$$
    where both the top horizontal and the left vertical maps are
    closed embeddings.  We will show that this diagram is a pushout
    square in the category of $k$-schemes.  In other words, for any
    $k$-scheme $Y$ forming a commutative diagram
    $$\xymatrix{Y & &\\
      & \spec k[s,t]/(s,t)^{2}\ar@{-->}[ul] & \spec 
      k[s,t]/(s^{2},t)\ar@(ul,r)[ull]\ar[l]\\
      & \spec k[s,t]/(s,t^{2})\ar@(l,d)[uul]\ar[u] & \spec 
      k,\ar[l]\ar[u]}$$
    there exists a unique $k$-morphism $\spec k[s,t]/(s,t)^{2}
    \rightarrow Y$ making the whole diagram commutative.  
    
    Since both the images of $\spec k[s,t]/(s^{2},t)$ and $\spec
    k[s,t]/(s,t^{2})$ in $Y$ are contained in some affine subscheme
    $Y_{0} \subseteq Y$, we can replace $Y$ with $Y_{0}$ and therefore
    assume $Y$ is affine.  Now the existence of the unique
    $k$-morphism $\spec k[s,t]/(s,t)^{2} \rightarrow Y$ is obvious
    because of the pullback square of $k$-algebras mentioned earlier.  So the
    scheme $\spec k[s,t]/(s,t)^{2}$ satisfies the universal property
    of a coproduct of $k$-schemes.
    
    Next, we apply the covariant functor $Z \times_{k} -$ to this
    pushout square of $k$-schemes, and one can check that we have a
    pushout square: 
    $$\xymatrix{Z \times_{k} \spec k[s,t]/(s,t)^{2} & Z \times_{k} 
      \spec k[s]/(s^{2})\ar[l]\\
      Z \times_{k} \spec k[t]/(t^{2})\ar[u] & Z \times_{k}
      \spec k\ar[l]\ar[u]}$$
    Finally, applying the contravariant and left exact functor
    $\Hom_{k}(-,X)$, we have a pullback square in the category of sets: 
    $$\xymatrix{\Hom_{k}(Z \times_{k} \spec k[s,t]/(s,t)^{2},X)\ar[d]\ar[r] 
      & \Hom_{k}(Z \times_{k} \spec k[s]/(s^{2}),X)\ar[d]\\
      \Hom_{k}(Z \times_{k} \spec k[t]/(t^{2}),X)\ar[r] &
      \Hom_{k}(Z,X)}$$
      
    On the other hand, we also have an obvious pullback square in the 
    category of {\bf k}\text{\bf-Schemes}:
    $$\xymatrix{\mathcal{J}_{1}(X) \times_{X} \mathcal{J}_{1}(X)
      \ar[d]\ar[r] & \mathcal{J}_{1}(X) \ar[d]\\
      \mathcal{J}_{1}(X) \ar[r]& X}$$
    We apply the covariant and left exact functor $\Hom_{k}(Z,-)$ and  
    obtain again a pullback square in the category of sets:
    $$\xymatrix{\Hom_{k}(Z,\mathcal{J}_{1}(X) \times_{X} \mathcal{J}_{1}(X)) 
      \ar[d]\ar[r] & \Hom_{k}(Z,\mathcal{J}_{1}(X))\ar[d]\\
      \Hom_{k}(Z,\mathcal{J}_{1}(X))\ar[r] & \Hom_{k}(Z,X)}$$
      
    By the uniqueness of fiber products, we have $$\Hom_{k}(Z \times \spec
    k[s,t]/(s,t)^{2},X) =\Hom_{k}(Z,\mathcal{J}_{1}(X) \times_{X}
    \mathcal{J}_{1}(X))$$ for all (affine) $k$-schemes $Z$. Therefore, 
    $\mathcal{W}_{1}(X) \cong \mathcal{J}_{1}(X) \times_{X} 
    \mathcal{J}_{1}(X)$.
\end{proof}

\section{Truncated wedge schemes of local complete
intersections}\label{IrreducibilityCriterion}

In the case of a locally complete intersection variety, we give an
irreducibility criterion for its truncated wedge schemes similar to
that of Musta\c{t}\v{a} for jet schemes \cite[Proposition
1.4]{MustataLCI}.

\begin{theorem}\label{LCI}
    Let $X$ be locally a complete intersection variety of dimension
    $n$.  Then the scheme $\mathcal{W}_m(X)$ is pure dimensional if
    and only if $$\dim
    \mathcal{W}_m(X)=\frac{1}{2}n(m+1)(m+2),$$ and in this case
    $\mathcal{W}_m(X)$ is locally a complete intersection.  Similarly,
    $\mathcal{W}_m(X)$ is irreducible if and only if $$\dim
    \pi_m^{-1}(X^{sing})<\frac{1}{2}n(m+1)(m+2),$$ where 
    $\pi_{m}:\mathcal{W}_{m}(X) \rightarrow X$ are the natural 
    projections.
\end{theorem}
\begin{proof}
    We have a decomposition
    \begin{equation}\label{Decomposition}
    \mathcal{W}_{m}(X)=\pi_{m}^{-1}(X^{sing}) \cup
    \overline{\pi_{m}^{-1}(X^{reg})}
    \end{equation}
    and in general $\overline{\pi_{m}^{-1}(X^{reg})}$ is an
    irreducible component of $\mathcal{W}_{m}(X)$ of dimension
    $\frac{1}{2}n(m+1)(m+2)$ by Proposition \ref{SmoothWedge}.  So the
    ``only if'' part of both assertions is obvious and holds without
    the local complete intersection hypothesis.

    Suppose that $\dim \mathcal{W}_{m}(X)=\frac{1}{2}n(m+1)(m+2)$.
    Working locally, we may assume that $X \subseteq \mathbb{A}^{N}$
    and $X$ is defined by $N-n$ equations.  Notice that each defining
    equation of $X$ gives rise to $\frac{1}{2}(m+1)(m+2)$ defining
    equations of $\mathcal{W}_{m}(X)$.  So $\mathcal{W}_{m}(X)
    \subseteq
    \mathcal{W}_{m}(\mathbb{A}^{N})=\mathbb{A}^{\frac{1}{2}N(m+1)(m+2)}$
    is defined by $\frac{1}{2}(N-n)(m+1)(m+2)$ equations.  Then by
    Krull's principal ideal theorem \cite[Theorem 10.2]{Eisenbud},
    every irreducible component of $\mathcal{W}_{m}(X)$ has dimension
    at least $\frac{1}{2}n(m+1)(m+2)$.  Thus $\dim
    \mathcal{W}_{m}(X)=\frac{1}{2}n(m+1)(m+2)$ implies that
    $\mathcal{W}_{m}(X)$ is pure dimensional and a local complete
    intersection.

    Now if $\dim \pi_{m}^{-1}(X^{sing})<\frac{1}{2}n(m+1)(m+2)$, the
    decomposition (\ref{Decomposition}) yields $\dim
    \mathcal{W}_{m}(X)=\dim
    \overline{\pi_{m}^{-1}(X^{reg})}=\frac{1}{2}n(m+1)(m+2)$.  So by
    what we just proved, the scheme $\mathcal{W}_{m}(X)$ is pure
    dimensional.  Therefore, our assumption on $\dim
    \pi_{m}^{-1}(X^{sing})$ tells us that $\pi_{m}^{-1}(X^{sing})$
    contributes no components to $\mathcal{W}_{m}(X)$.  Hence,
    $\pi_{m}^{-1}(X^{sing})$ is contained in the closure of
    $\pi_{m}^{-1}(X^{reg})$, and $\mathcal{W}_{m}(X)$ is irreducible.
\end{proof}

\section{Truncated wedge schemes of monomial schemes}

In this section, we analyze the scheme structure of truncated wedge
schemes of monomial schemes.  Like jet schemes, truncated wedge
schemes of monomial schemes are not themselves monomial but our
calculations will show that their reduced subschemes are. The 
following theorem is analogous to that of Goward and Smith for jet 
schemes \cite[Theorem 3.1]{GowardSmith}.

\begin{theorem}\label{MonomialRadical}
    Let $X \subseteq \mathbb{A}^{n}$ be a monomial scheme in
    coordinates $x_{1},\ldots,x_{n}$.  Then $\sqrt{W_{m}(X)}$ is a
    square-free monomial ideal in the coordinates
    $x_{k}^{(i_{k},j_{k})}$ where $1 \leq k \leq n$ and $0 \leq
    i_{k}+j_{k} \leq m$.  The generators of $\sqrt{W_{m}(X)}$ can be
    described as follows: for each minimal monomial generator of the
    defining ideal of $X$, say $x_{1}^{a_{1}}\cdots x_{r}^{a_{r}}$
    after relabeling, the monomials
    $$\sqrt{x_{1}^{(i_{1},j_{1})}x_{1}^{(i_{2},j_{2})}\cdots
    x_{1}^{(i_{a_{1}},j_{a_{1}})}x_{2}^{(i_{a_{1}+1},j_{a_{1}+1})}
    \cdots x_{2}^{(i_{a_{1}+a_{2}},j_{a_{1}+a_{2}})}\cdots
    x_{r}^{(i_{a_{1}+\ldots+a_{r}},j_{a_{1}+\ldots+a_{r}})}}$$
    with $\sum (i_{k}+j_{k}) \leq m$ are monomial generators of
    $\sqrt{W_{m}(X)}$.  The collection of all such monomials as we
    range through the minimal monomial generators of the defining
    ideal of $X$ is a generating set for $\sqrt{W_{m}(X)}$.
\end{theorem}

The following lemma taken from \cite[Lemma 2.1]{GowardSmith} allows us
to reduce the proof of Theorem \ref{MonomialRadical} to the
hypersurface case.

\begin{lemma}
    If $I$ and $J$ are monomial ideals in a polynomial ring, then
    $\sqrt{I+J}=\sqrt{I}+\sqrt{J}$.
\end{lemma}

\subsection{\textit{Monomial hypersurface case}}

\begin{theorem}\label{GeneralMonomialHypersurface}
    Let $X=\spec k[x_{1},\ldots,x_{n}]/(x_{1}^{a_{1}}\cdots
    x_{r}^{a_{r}})$ be a monomial hypersurface.  Then
    \begin{enumerate}
	\item The minimal primes of $W_{m}(X)$ are
	precisely the minimal members of the set of prime ideals
	\begin{equation}\label{MinimalPrimes}
	    (x_{k}^{(i_{k},j_{k})}:0 \leq i_{k}+j_{k}<t_{k},1 \leq k
	    \leq r)
	\end{equation}
	where $0 \leq t_{k} \leq m+1$ and $\sum a_{k}t_{k}\geq m+1$.
	(Here, we use the convention that the value $t_{k}=0$ means
	that the variable $x_{k}$ does not appear at all.)  \label{part1}
	\item The radical $\sqrt{W_{m}(X)}$ is the
	monomial ideal generated by the monomials
	$$\sqrt{x_{1}^{(i_{1},j_{1})}x_{1}^{(i_{2},j_{2})}\cdots
	x_{1}^{(i_{a_{1}},j_{a_{1}})}x_{2}^{(i_{a_{1}+1},j_{a_{1}+1})}
	\cdots x_{2}^{(i_{a_{1}+a_{2}},j_{a_{1}+a_{2}})}\cdots
	x_{r}^{(i_{a_{1}+\ldots+a_{r}},j_{a_{1}+\ldots+a_{r}})}}$$
	where $\sum (i_{k}+j_{k}) \leq m$. \label{part2}
    \end{enumerate}	
\end{theorem}

Before we prove this result, we first need to describe the generators
of $W_{m}(X)$.  Recall that the polynomials defining the truncated
wedge scheme $\mathcal{W}_{m}(X)$ are the coefficients of $s^{i}t^{j}$
in the product
$$\prod_{k=1}^{r}\left(x_{k}^{(0,0)}+x_{k}^{(1,0)}s+x_{k}^{(0,1)}t
+\ldots+x_{k}^{(m,0)}s^{m}+x_{k}^{(m-1,1)}s^{m}t+\ldots
+x_{k}^{(0,m)}t^{m}\right)^{a_{k}}.$$ Therefore, $W_{m}(X)$ has generators 
of the form 
\begin{equation}\label{g_{ij}}
    g_{ij}=\sum x_{1}^{(i_{1},j_{1})}\cdots
    x_{1}^{(i_{a_{1}},j_{a_{1}})}x_{2}^{(i_{a_{1}+1},j_{a_{1}+1})}
    \cdots x_{2}^{(i_{a_{1}+a_{2}},j_{a_{1}+a_{2}})}\cdots
    x_{r}^{(i_{a_{1}+\ldots+a_{r}},j_{a_{1}+\ldots+a_{r}})}
\end{equation}    
where $0 \leq i+j \leq m$, $0 \leq i_{k},j_{k} \leq m$, $\sum i_{k}=i$
and $\sum j_{k}=j$.

\begin{proof}[Proof of Theorem \ref{GeneralMonomialHypersurface}]
    We will first show that part (\ref{part2}) follows from part 
    (\ref{part1}).
    Recall that the radical of any ideal is equal to the intersection
    of its minimal primes.  So by part (\ref{part1}), it is clear that
    $\sqrt{W_{m}(X)}$ is a monomial ideal.  Next, observe that the
    monomials
    \begin{equation}\label{Monomials}
	x_{1}^{(i_{1},j_{1})}x_{1}^{(i_{2},j_{2})}\cdots
	x_{1}^{(i_{a_{1}},j_{a_{1}})}x_{2}^{(i_{a_{1}+1},j_{a_{1}+1})}
	\cdots x_{2}^{(i_{a_{1}+a_{2}},j_{a_{1}+a_{2}})}\cdots
	x_{r}^{(i_{a_{1}+\ldots+a_{r}},j_{a_{1}+\ldots+a_{r}})}
    \end{equation}
    with $\sum (i_{k}+j_{k}) \leq m$ are precisely the terms of the
    polynomials $g_{ij}$ as $i+j$ ranges from $0$ to $m$.  Since every
    monomial ideal must contain all terms of its generators, it
    follows that the monomials in (\ref{Monomials}) and therefore the
    monomials listed in part (\ref{part2}) are all contained in
    $\sqrt{W_{m}(X)}$.  On the other hand, since these monomials are
    all square-free, they generate a radical ideal containing
    $W_{m}(X)$.  Thus, this is the smallest radical ideal containing
    $W_{m}(X)$, and hence must be $\sqrt{W_{m}(X)}$ exactly.
     
    To prove part (\ref{part1}), we induce on $m$.  The ideal $W_{0}(X)$ is
    $\left((x_{1}^{(0,0)})^{a_{1}}\cdots
    (x_{r}^{(0,0)})^{a_{r}}\right)$ and its minimal primes are
    obviously $(x_{i}^{(0,0)})$ for $1 \leq i \leq r$.  For the
    inductive step, let $Q$ be a prime containing $W_{m}(X)$.  Since
    $W_{m-1}(X) \subset W_{m}(X)$, the prime $Q$ contains one of the
    minimal primes $P$ of $W_{m-1}(X)$.  So by induction,
    $$P=(x_{k}^{(i_{k},j_{k})}:0 \leq i_{k}+j_{k}<t_{k},1 \leq k \leq
    r)$$ for some $0 \leq t_{k} \leq m$ and $\sum a_{k}t_{k}\geq m$.
    Fix the indices $t_{1},\ldots,t_{r}$ corresponding to $P$.  We
    want to show $Q$ contains a ``full layer'' of variables; that is,
    a set of elements $\{x_{n}^{(i_{n},t_{n}-i_{n})}:0 \leq i_{n} \leq
    t_{n}\}$ for some $n \in \{1,\ldots,r\}$.

    Consider the set $\mathcal{S}$ of primes lying between $P$ and
    $Q$, and containing several ``initial partial layers''.  In other
    words, consider primes contained in $Q$ and of the form
    $$P+(x_{k}^{(i_{k},t_{k}-i_{k})}:0 \leq i_{k}<c_{k}, 1 \leq k \leq
    r)$$ for some $0 \leq c_{k} \leq t_{k}$.
    
    Note that the set $\mathcal{S}$ is nonempty because $P$ is a
    member (take $c_{k}$ to be zero for all $k$).  Now choose a
    maximal element $\mathcal{P}$ in the set $\mathcal{S}$
    ($\mathcal{P}$ exists by the Noetherian property), and consider
    the monomial $$\alpha=(x_{1}^{(c_{1},t_{1}-c_{1})})^{a_{1}}\cdots
    (x_{r}^{(c_{r},t_{r}-c_{r})})^{a_{r}}.$$ Clearly,
    $\alpha \notin \mathcal{P}$ and $\alpha$ is a monomial term of the
    generator $g_{c,m-c}$ where $c=\sum_{k=1}^{r} a_{k}c_{k}$.
    \begin{claim}\label{g_{c,m-c}}
	All monomial terms of $g_{c,m-c}$ except $\alpha$ belong to
	the prime $\mathcal{P}$.
    \end{claim}
    \begin{proof}
	Let $\beta$ be a monomial term of $g_{c,m-c}$.  Then
	$$\beta=x_{1}^{(i_{1},j_{1})}\cdots
	x_{1}^{(i_{a_{1}},j_{a_{1}})}x_{2}^{(i_{a_{1}+1},j_{a_{1}+1})}
	\cdots x_{2}^{(i_{a_{1}+a_{2}},j_{a_{1}+a_{2}})}\cdots
	x_{r}^{(i_{a_{1}+\ldots+a_{r}},j_{a_{1}+\ldots+a_{r}})}$$
	where $\sum i_{l}=c$ and $\sum j_{l}=m-c$.  Suppose $\beta
	\notin \mathcal{P}$, and so in particular, $\beta \notin P$.
	Then the description of $P$ implies that every sum of
	superscripts $i_{l}+j_{l}$ attached to each $x_{k}$ is at
	least $t_{k}$ for all $k=1,\ldots,r$.  So $$m=\sum (i_{l}+j_{l}) \geq
	\sum a_{k}t_{k} \geq m$$ implies that $i_{l}+j_{l}=t_{k}$ for
	all $k$ and $l$.  If $i_{l} \geq c_{k}$ for all $k$ and $l$,
	then $$c=\sum i_{l} \geq \sum a_{k}c_{k}=c$$
	implies that $i_{l}=c_{k}$ for all $k,l$.  So
	$$\beta=(x_{1}^{(c_{1},t_{1}-c_{1})})^{a_{1}}\cdots
	(x_{r}^{(c_{r},t_{r}-c_{r})})^{a_{r}}=\alpha.$$ Thus, the 
	only monomial term of $g_{c,m-c}$ not in $\mathcal{P}$ is 
	$\alpha$. This completes the proof of Claim \ref{g_{c,m-c}}.
    \renewcommand{\qed}{}
    \end{proof}

    Since $g_{c,m-c} \in W_{m}(X) \subset Q$, Claim \ref{g_{c,m-c}}
    gives $\alpha \in Q$ and so $x_{n}^{(c_{n},t_{n}-c_{n})} \in Q$
    for some $n \in \{1,\ldots,r\}$ by primality of $Q$.  Fix such an
    $n$.  If $c_{n}<t_{n}$, then
    $\mathcal{P}+(x_{n}^{(c_{n},t_{n}-c_{n})})$ is an element of the
    set $\mathcal{S}$, violating maximality of $\mathcal{P}$.
    Therefore, $c_{n}=t_{n}$, which means $Q$ contains the ideal
    $$P':=P+(x_{n}^{(i_{n},t_{n}-i_{n})}:0 \leq i_{n} \leq t_{n}),$$ a
    complete layer, as desired.

    We are left to check that $W_{m}(X)$ is contained in the ideal
    $P'$.  Since $W_{m-1}(X) \subset P$, the only generators of
    $W_{m}(X)$ not already in $P$ are $g_{a,m-a}$ for $0 \leq a \leq
    m$.  By a similar argument as in the proof of Claim
    \ref{g_{c,m-c}}, we see that a term of $g_{a,m-a}$ not in
    $P$ must have every sum of superscripts $i_{l}+j_{l}$ attached to
    each $x_{k}$ exactly $t_{k}$.  Obviously, this term is in the
    ideal $(x_{n}^{(i_{n},t_{n}-i_{n})}:0 \leq i_{n} \leq t_{n})$.
    Since every prime containing $W_{m}(X)$ contains a prime of the
    form $P'$, the minimal primes of $W_{m}(X)$ are precisely those
    described in (\ref{MinimalPrimes}).
\end{proof}

In the special case when $X$ is a reduced monomial hypersurface, we 
get a slightly sharper result:
\begin{corollary}\label{ReducedMonomialHypersurface}
    Let $X=\spec k[x_{1},\ldots,x_{n}]/(x_{1}\cdots x_{r})$ be a
    reduced hypersurface monomial scheme.  Then the minimal primes of
    $W_{m}(X)$ are exactly the primes of the form
    $$P=(x_{k}^{(i_{k},j_{k})}:0 \leq i_{k}+j_{k}<t_{k},1 \leq k \leq
    r)$$ where $0 \leq t_{k} \leq m+1$ and $\sum t_{k}=m+1$.  (Again,
    $t_{k}=0$ means the variable $x_{k}$ does not appear.)
\end{corollary}

\begin{rmk}\label{NotPureDimensional}
    This result tells us that truncated wedge schemes of a monomial 
    scheme are not in general pure dimensional. For example, take 
    $r=3$ and $m=2$ in Corollary \ref{ReducedMonomialHypersurface}. 
    Then $W_{2}(X)$ has a minimal prime of height $6$ (for example, 
    take $t_{1}=3$ and $t_{2}=t_{3}=0$), a minimal prime of height 
    $4$ (take $t_{1}=2$, $t_{2}=1$, $t_{3}=0$), and a minimal prime 
    of height $3$ (take $t_{1}=t_{2}=t_{3}=1$).
\end{rmk}    

\subsection{\textit{Alternate point of view}}\label{GeometricApproach}

In this section, we present another approach towards Theorem
\ref{GeneralMonomialHypersurface}.  We will examine
$\mathcal{W}_{m}(X)$ first as a set; that is, we will look at the
reduced subscheme of $\mathcal{W}_{m}(X)$ and then describe the
irreducible components of $\mathcal{W}_{m}(X)$.  For $X=\spec
k[x_{1},\ldots,x_{n}]/(x_{1}^{a_{1}}\cdots x_{r}^{a_{r}})$, a point in
$\mathcal{W}_{m}(X)$ corresponds to an $N$-tuple
$(x_{l}^{(i_{l},j_{l})}) \in k^{N}$, $N=\frac{1}{2}r(m+1)(m+2)$,
satisfying $$\prod_{l=1}^{r} (\sum x_{l}^{(i_{l},j_{l})}s^{i_{l}}
t^{j_{l}})^{a_{l}} \in (s,t)^{m+1}.$$ 

We write $\ord(f)$ for the smallest positive integer $n$ such that $f
\in (s,t)^{n}$, and note that $\ord(fg)=\ord(f)+\ord(g)$.  Then
writing $p_{l}$ for $\sum x_{l}^{(i_{l},j_{l})}s^{i_{l}}t^{j_{l}}$,
the condition that $$\prod_{l=1}^{r} p_{l}^{a_{l}} \in (s,t)^{m+1}$$
means $\ord(\prod p_{l}^{a_{l}}) \geq m+1$, or equivalently, $\sum
a_{l}\ord(p_{l}) \geq m+1$.  Now let
\begin{align*}
    t_{l} &:= \ord(p_{l})\\
    &= \min \{c:x_{l}^{(i_{l},j_{l})}=0 \text{ for all } i_{l}+j_{l}<c\}.
\end{align*}
Then a point in $\mathcal{W}_{m}(X)$ corresponds to an $N$-tuple
$(x_{l}^{(i_{l},j_{l})}) \in k^{N}$ satisfying
$x_{l}^{(i_{l},j_{l})}=0$ for $1 \leq l \leq r$, $0 \leq
i_{l}+j_{l}<t_{l}$, and $\sum a_{l}t_{l} \geq m+1$.  Thus the
vanishing set of the generators of $W_{m}(X)$ is
\begin{align*}
    \mathbb{V}(g_{ij}) &= \bigcup_{\sum a_{l}t_{l} \geq m+1}
    \mathbb{V}(x_{l}^{(i_{l},j_{l})}:0 \leq i_{l}+j_{l}<t_{l},1 \leq l
    \leq r)\\
    &= \mathbb{V}\left(\bigcap_{\sum a_{l}t_{l} \geq m+1}
    (x_{l}^{(i_{l},j_{l})}:0 \leq i_{l}+j_{l}<t_{l},1 \leq l \leq
    r)\right).
\end{align*}

Since each ideal $(x_{l}^{(i_{l},j_{l})}:0 \leq i_{l}+j_{l}<t_{l},1
\leq l \leq r)$ is prime, their intersection is a radical ideal.
Therefore, the radical
\begin{equation}\label{Radical}
    \sqrt{W_{m}(X)}=\bigcap_{\sum a_{l}t_{l} \geq m+1}
    (x_{l}^{(i_{l},j_{l})}:0 \leq i_{l}+j_{l}<t_{l},1 \leq l \leq r)
\end{equation}
is a monomial ideal.  Since the right hand side of (\ref{Radical}) is
a primary decomposition of $\sqrt{W_{m}(X)}$, the minimal primes of
$W_{m}(X)$ are precisely the minimal members of the collection of
primes $$\left\{(x_{k}^{(i_{k},j_{k})}: 0 \leq i_{k}+j_{k}<t_{k}, \sum
a_{k}t_{k} \geq m+1)\right\}.$$

\subsection{\textit{Multiplicity}}\label{WedgeMultiplicity}

Although the defining ideal and its minimal primes of a truncated
wedge scheme of a reduced monomial hypersurface $X$ are very similar
to those of a jet scheme, there are important differences.  From
Remark \ref{NotPureDimensional}, we see that the truncated wedge
schemes of $X$ are not pure dimensional, while the jet schemes are
\cite[Theorem 2.2]{GowardSmith}.  Multiplicities behave quite
differently too.  It was shown in \cite{Yuen} that the multiplicity of
the jet schemes $\mathcal{J}_{m}(X)$ along a component is a
multinomial coefficient, depending on the component.  But Macaulay
calculations suggest that the multiplicity of $\mathcal{W}_{m}(X)$
along any irreducible component is always one.

\begin{conj}\label{WedgeMultConj}
    If $X=\spec k[x_{1},\ldots,x_{n}]/(x_{1}\cdots x_{r})$ is a
    reduced monomial hypersurface, then the multiplicity of
    $\mathcal{W}_{m}(X)$ along any irreducible component is one.
\end{conj}

\begin{rmk}\label{EmbeddedComponent}
    This conjecture does not imply the schemes $W_{m}(X)$ are reduced
    because they may have embedded components.  For example, a
    Macaulay calculation shows that for the scheme $X=\spec
    k[x,y]/(xy)$, the ideal $W_{1}(X)$ has three minimal primes
    $$(x_{00},y_{00}), \hspace{0.3cm} (x_{00},x_{01},x_{10}),
    \hspace{0.3cm}(y_{00},y_{01},y_{10})$$
    and one non-minimal associated prime
    $(x_{00},y_{00},x_{01}y_{10}-x_{10}y_{01})$.
\end{rmk}

We will prove Conjecture \ref{WedgeMultConj} in the case when $r=2$ or
$3$.  We will show that a certain subset of the generators of
$W_{m}(X)$ generate the maximal ideal of the local ring along the
corresponding irreducible component.

\begin{theorem}\label{WedgeMultI}
    For $X=\spec k[x,y,z_{1},\ldots,z_{n}]/(xy)$, the multiplicity of
    $\mathcal{W}_{m}(X)$ along any irreducible component is one.
\end{theorem}
\begin{proof}
    Fix a minimal prime $P$ of $W_{m}(X)$, that is,
    $$P=(x_{ij},y_{kl}:0 \leq i+j<t_{1},0 \leq k+l<t_{2})$$ for some
    $0 \leq t_{1},t_{2} \leq m+1$ such that $t_{1}+t_{2}=m+1$.  For
    each $x_{ij} \in P$, pick the generator $g_{i,j+t_{2}}$, and for
    each $y_{kl} \in P$, pick the generator $g_{k+t_{1},l}$.  These
    are generators of $W_{m}(X)$ because $0 \leq i+j<t_{1}$, $0
    \leq k+l<t_{2}$ and $t_{1}+t_{2}=m+1$ imply that $0 \leq
    i+(j+t_{2}),(k+t_{1})+l \leq m$.  Also, note that $g_{i,j+t_{2}}
    \neq g_{k+t_{1},l}$ since $i<t_{1}$.
    
    We want to show that this subset of generators form a regular 
    system of parameters in $R_{P}$ where $R$ is the coordinate ring 
    of $\mathcal{W}_{m}(\spec k[x,y,z_{1},\ldots,z_{n}])$. Indeed, we 
    will show that the images of these elements in 
    $PR_{P}/P^{2}R_{P}$ form a basis for the cotangent space 
    $P/P^{2}$ of the regular local ring $R_{P}$.
    
    Let $M$ be the matrix whose rows express the elements
    $\overline{g_{ab}}$ as $R_{P}/PR_{P}$-linear combination of the
    basis $\overline{x_{ij}}$, $\overline{y_{kl}}$ of the vector space
    $P/P^{2}$.  Then from the expression (\ref{g_{ij}}), we can check
    that the matrix $M$ has the form
    
    $$\begin{array}{cc@{}cccccccc@{}c}
      & & \overline{x_{00}} & \overline{x_{10}} & \cdots &
      \overline{x_{0,t_{1}-1}} & \overline{y_{00}} & \overline{y_{10}}
      & \cdots & \overline{y_{0,t_{2}-1}}\\
      \begin{matrix} \overline{g_{0,t_{2}}}\\ \overline{g_{1,t_{2}}}\\
	  \vdots \\ \overline{g_{0,t_{1}+t_{2}-1}}\\ \overline{g_{t_{1},0}}\\
	  \overline{g_{t_{1}+1,0}}\\ \vdots \\\overline{g_{t_{1},t_{2}-1}}
          \end{matrix} &
      \left(\vphantom{\begin{matrix} \overline{g_{0,t_{2}}}\\ 
          \overline{g_{1,t_{2}}}\\ \vdots \\ 
	  \overline{g_{0,t_{1}+t_{2}-1}}\\ \overline{g_{t_{1},0}}\\ 
	  \overline{g_{t_{1}+1,0}}\\ \vdots \\ \overline{g_{t_{1},t_{2}-1}} 
          \end{matrix}}\right. &
      \begin{matrix} \overline{y_{0,t_{2}}}\\ \star \\ \vdots \\ \star \\
          \star \\ \star \\ \vdots \\ \star \end{matrix} &
      \begin{matrix} \star \\ \overline{y_{0,t_{2}}}\\ \vdots \\ \star \\
          \star \\ \star \\ \vdots \\ \star \end{matrix} &
      \begin{matrix} \cdots \\ \cdots \\ \ddots \\ \cdots \\
          \cdots \\ \cdots \\ \ddots \\ \cdots \end{matrix} &
      \begin{matrix} \star \\ \star \\ \vdots \\ \overline{y_{0,t_{2}}}\\
          \star \\ \star \\ \vdots \\ \star \end{matrix} &
      \begin{matrix} \star \\ \star \\ \vdots \\ \star \\
          \overline{x_{t_{1},0}}\\ \star \\ \vdots \\ \star \end{matrix} &
      \begin{matrix} \star \\ \star \\ \vdots \\ \star \\
          \star \\ \overline{x_{t_{1},0}}\\ \vdots \\ \star \end{matrix} &
          \begin{matrix} \cdots \\ \cdots \\ \ddots \\ \cdots \\
          \cdots \\ \cdots \\ \ddots \\ \cdots \end{matrix} &
      \begin{matrix} \star \\ \star \\ \vdots \\ \star \\
          \star \\ \star \\ \vdots \\ \overline{x_{t_{1},0}} \end{matrix} &
      \left.\vphantom{\begin{matrix} \overline{g_{0,t_{2}}}\\ 
          \overline{g_{1,t_{2}}}\\ \vdots \\ 
	  \overline{g_{0,t_{1}+t_{2}-1}}\\ \overline{g_{t_{1},0}}\\ 
	  \overline{g_{t_{1}+1,0}}\\ \vdots \\ \overline{g_{t_{1},t_{2}-1}} 
          \end{matrix}}\right)
      \end{array}$$
    where $\star$'s are of the form $\overline{x_{ab}}$ where $(a,b)
    \neq (t_{1},0)$ or $\overline{y_{ab}}$ where $(a,b) \neq
    (0,t_{2})$.

    Now to show that the matrix $M$ is invertible, note that its
    determinant is a polynomial in the variables $\overline{x_{ab}}$
    and $\overline{y_{cd}}$ for $a+b \geq t_{1}$ and $c+d \geq t_{2}$,
    interpreted as an element in the field $$R_{P}/PR_{P} \cong
    k(x_{ab}, y_{cd}:t_{1} \leq a+b \leq m, t_{2} \leq c+d \leq m).$$
    This polynomial is not the zero polynomial because plugging in the
    values $\overline{x_{t_{1},0}}=\overline{y_{0,t_{2}}}=1$ and all
    other $\overline{x_{ab}}$ and $\overline{y_{cd}}$ to be zero, we
    get the nonzero value $1$.  So the matrix $M$ is invertible, which
    means the elements
    $\overline{g_{0,t_{2}}},\ldots,\overline{g_{t_{1},t_{2}-1}}$ form
    a basis for the $R_{P}/PR_{P}$-vector space $P/P^{2}$.  Therefore,
    the generators $g_{0,t_{2}},\ldots,g_{t_{1},t_{2}-1}$ form a
    regular system of parameters in $R_{P}$.  Thus, $W_{m}(X)R_{P}
    \supseteq (g_{0,t_{2}},\ldots,g_{t_{1},t_{2}-1})R_{P}=PR_{P}$, and
    so $\ell(R_{P}/W_{m}(X))=1$.
\end{proof}

\begin{theorem}\label{WedgeMultII}
    For $X=\spec k[x,y,z,w_{1},\ldots,w_{n}]/(xyz)$, the multiplicity of
    $\mathcal{W}_{m}(X)$ along any irreducible component is one.
\end{theorem}
\begin{proof}
    Consider a minimal prime $$P=(x_{ij},y_{kl},z_{pq}:0 \leq
    i+j<t_{1},0 \leq k+l<t_{2},0 \leq p+q<t_{3})$$ for some $0 \leq t_{r}
    \leq m+1$ and $\sum t_{r}=m+1$.  Without loss of generality, we
    may assume $t_{1} \leq t_{2} \leq t_{3}$.  For each $x_{ij} \in
    P$, pick the generator $g_{i,j+t_{2}+t_{3}}$.  For each $y_{kl}
    \in P$, pick the generator $g_{k+t_{1}+t_{3},l}$.  For each
    $z_{pq} \in P$, pick the generator $g_{p+t_{1},q+t_{2}}$.  Notice
    that these are generators of $W_{m}(X)$, since $0 \leq i+j<t_{1}$,
    $0 \leq k+l<t_{2}$, $0 \leq p+q<t_{3}$ and
    $\sum t_{r}=m+1$ imply that $0 \leq
    i+(j+t_{2}+t_{3}),(k+t_{1}+t_{3})+l,(p+t_{1})+(q+t_{2}) \leq m$.
    Also, they are distinct because $i<t_{1}$ and $l<t_{2}$.

    We want to show that the images of these generators in
    $PR_{P}/P^{2}R_{P}$ form a basis for the cotangent space $P/P^{2}$
    of the regular local ring $R_{P}$, where $R$ is the coordinate
    ring of the scheme $\mathcal{W}_{m}(\spec
    k[x,y,z,w_{1},\ldots,w_{n}])$.  As in the previous proof, let $M$
    be the matrix whose rows express the elements $\overline{g_{ab}}$
    as $R_{P}/PR_{P}$-linear combination of the basis
    $\overline{x_{ij}}$, $\overline{y_{kl}}$ and $\overline{z_{pq}}$
    of the vector space $P/P^{2}$.  We will show that $M$ is
    invertible.

    From the expression (\ref{g_{ij}}), one can check that the matrix
    $M$ has the form 
    $$\begin{array}{cc@{}ccccccccc@{}c} 
          & & \hphantom{x_{ij}} & \overline{x_{ij}} & \hphantom{x_{ij}} &
	  \hphantom{y_{kl}} & \overline{y_{kl}} & \hphantom{y_{kl}} &
	  \hphantom{z_{pq}} & \overline{z_{pq}} & \hphantom{z_{pq}}\\
          \begin{matrix} \\ \overline{g_{i,j+t_{2}+t_{3}}}\\ \\ \\
	      \overline{g_{k+t_{1}+t_{3},l}}\\ \\ \\
	      \overline{g_{p+t_{1},q+t_{2}}}\\ \\ \end{matrix} &
	  \left(\vphantom{\begin{matrix} \\ g_{i,j+t_{2}+t_{3}}\\ \\
	      \\ g_{k+t_{1}+t_{3},l}\\ \\ \\
	      g_{p+t_{1},q+t_{2}}\\ \\ \end{matrix}}\right.  &
	  \begin{matrix} h_{1}+\star \\
	      \vdots \\ \star \\ \star \\ \vdots \\ \star \\ \heartsuit
	      \\ \vdots \\ \heartsuit \end{matrix} &
	  \begin{matrix} \cdots \\ \ddots \\ \cdots \\ \cdots \\
	      \ddots \\ \cdots \\ \cdots \\ \ddots \\ \dots \end{matrix} &
	  \begin{matrix} \heartsuit \\ \vdots \\
	      h_{1}+\star \\ \star \\ \vdots \\
	      \star \\ \heartsuit \\ \vdots \\ \heartsuit \end{matrix} &
	  \begin{matrix} \star \\ \vdots \\ \star \\ h_{2}+\star
	       \\ \vdots \\ \heartsuit \\
	      \heartsuit \\ \vdots \\ \heartsuit \end{matrix} &
	  \begin{matrix} \cdots \\ \ddots \\ \cdots \\ \cdots \\
	      \ddots \\ \cdots \\ \cdots \\ \ddots \\ \dots \end{matrix} &
	  \begin{matrix} \star \\ \vdots \\ \star \\ \star \\ \vdots 
	      \\ h_{2}+\star \\
	      \heartsuit \\ \vdots \\ \heartsuit \end{matrix} & 
	  \begin{matrix} \star \\ \vdots \\ \star \\ \star \\ \vdots 
	      \\ \star \\ h_{3}+\star \\
	      \vdots \\ \star \end{matrix} & 
	  \begin{matrix} \cdots \\ \ddots \\ \cdots \\ \cdots \\
	      \ddots \\ \cdots \\ \cdots \\ \ddots \\ \dots \end{matrix} &
	  \begin{matrix} \star \\ \vdots \\ \star \\ \star \\ \vdots 
	      \\ \star \\ \star \\ \vdots \\ 
	      h_{3}+\star \end{matrix} & 
	  \left.\vphantom{\begin{matrix} \\ g_{i,j+t_{2}+t_{3}}\\ \\ \\
	      g_{k+t_{1}+t_{3},l}\\ \\ \\
	      g_{p+t_{1},q+t_{2}}\\ \\ \end{matrix}}\right) 
      \end{array}$$
    where $h_{1}=\overline{y_{0,t_{2}}z_{0,t_{3}}}$,
    $h_{2}=\overline{x_{t_{1},0}z_{t_{3},0}}$, and
    $h_{3}=\overline{x_{t_{1},0}y_{0,t_{2}}}$.  The symbol $\star$
    represents a polynomial in the variables $x_{ab}$ where $(a,b)
    \neq (t_{1},0)$, $y_{cd}$ where $(c,d) \neq (0,t_{2})$, and
    $z_{ef}$ where $(e,f) \neq (0,t_{3})$, $(t_{3},0)$.  The symbol
    $\heartsuit$ represents an arbitrary polynomial in $R_{P}/PR_{P}$.
    Interpreting the entries of this matrix as elements of the field
    $$R_{P}/PR_{P} \cong k(x_{ab}, y_{cd}, z_{ef}:t_{1} \leq a+b \leq
    m, t_{2} \leq c+d \leq m, t_{3} \leq e+f \leq m),$$ we see that to
    check the determinant of $M$ is nonzero, it suffices to check it
    has nonzero value for some choice of values for $x_{ab}$, $y_{cd}$
    and $z_{ef}$.  Now setting
    $\overline{x_{t_{1},0}}=\overline{y_{0,t_{2}}}=\overline{z_{0,t_{3}}}
    =\overline{z_{t_{3},0}}=1$ and all the others to be zero, we get
    $$\begin{array}{cc@{}ccccccccc@{}c} 
          & & \hphantom{x_{ij}} & \overline{x_{ij}} & \hphantom{x_{ij}} &
	  \hphantom{y_{kl}} & \overline{y_{kl}} & \hphantom{y_{kl}} &
	  \hphantom{z_{pq}} & \overline{z_{pq}} & \hphantom{z_{pq}}\\
          \begin{matrix} \\ \overline{g_{i,j+t_{2}+t_{3}}}\\ \\ \\
	      \overline{g_{k+t_{1}+t_{3},l}}\\ \\ \\
	      \overline{g_{p+t_{1},q+t_{2}}}\\ \\ \end{matrix} &
	  \left(\vphantom{\begin{matrix} \\ g_{i,j+t_{2}+t_{3}}\\ \\
	      \\ g_{k+t_{1}+t_{3},l}\\ \\ \\
	      g_{p+t_{1},q+t_{2}}\\ \\ \end{matrix}}\right.  &
	  \begin{matrix} 1\\ \vdots \\ 0\\ 0\\ \vdots \\ 0\\ \heartsuit \\
	      \vdots \\ \heartsuit \end{matrix} & 
	  \begin{matrix} \cdots \\ \ddots \\ \cdots \\ \cdots \\
	      \ddots \\ \cdots \\ \cdots \\ \ddots \\ \dots \end{matrix} &
	  \begin{matrix} \heartsuit \\ \vdots \\ 1\\ 0\\ \vdots \\ 0\\
	      \heartsuit \\ \vdots \\ \heartsuit \end{matrix} & 
	  \begin{matrix} 0\\ \vdots \\ 0\\ 1\\ \vdots \\ \heartsuit \\ \heartsuit \\
	      \vdots \\ \heartsuit \end{matrix} & 
	  \begin{matrix} \cdots \\ \ddots \\ \cdots \\ \cdots \\
	      \ddots \\ \cdots \\ \cdots \\ \ddots \\ \dots \end{matrix} &
	  \begin{matrix} 0\\ \vdots \\ 0\\ 0\\ \vdots \\ 1\\
	      \heartsuit \\ \vdots \\ \heartsuit \end{matrix} & 
	  \begin{matrix} 0\\ \vdots \\ 0\\ 0\\ \vdots \\ 0\\ 1\\
	      \vdots \\ 0 \end{matrix} & 
	  \begin{matrix} \cdots \\ \ddots \\ \cdots \\ \cdots \\
	      \ddots \\ \cdots \\ \cdots \\ \ddots \\ \dots \end{matrix} &
	  \begin{matrix} 0\\ \vdots \\ 0\\ 0\\ \vdots \\ 0\\ 0\\
	  \vdots \\ 1 \end{matrix} & 
	  \left.\vphantom{\begin{matrix} \\ g_{i,j+t_{2}+t_{3}}\\ \\ \\
	      g_{k+t_{1}+t_{3},l}\\ \\ \\
	      g_{p+t_{1},q+t_{2}}\\ \\ \end{matrix}}\right) 
      \end{array}$$
    where the $\heartsuit$'s represent arbitrary elements of $k$.

    By reversing the order of $\overline{x_{ij}}$ and therefore the
    order of $\overline{g_{i,j+t_{2}+t_{3}}}$, we obtain a
    lower-triangular matrix with $1$ on the diagonal.  So the above
    matrix has determinant $1$, and therefore the matrix $M$ is
    invertible.  Thus, the set
    $$\{g_{i,j+t_{2}+t_{3}},g_{k+t_{1}+t_{3},l},g_{p+t_{1},q+t_{2}}:0 
    \leq i+j<t_{1},0 \leq k+l<t_{2},0 \leq p+q<t_{3}\}$$
    is a regular system of parameters in $R_{P}$, and hence $\ell
    (R_{P}/W_{m}(X))=1$.
\end{proof}

\begin{rmk}
    In the case when $r=4$, Conjecture \ref{WedgeMultConj} was checked
    up to $m=5$ by Macaulay 2.
\end{rmk}

\bibliographystyle{plain}
\bibliography{biblio}   % Use the BibTeX file ``biblio.bib''.

\end{document}